\documentclass{amsart}
\usepackage{amsfonts}
\usepackage{mathrsfs}  % [00pt]changes size of font
 \usepackage{amsmath}    % math enhancements latex2e (replaces amstex)
 \usepackage{amsthm}     % proclaims theoremstyle/proof environments
 \usepackage{amscd}      % commutative Diagram environment
 \usepackage{amssymb}    % AMSFonts and symbols
 \usepackage{eucal}      % Euler Cal/Script Fonts
 \usepackage{latexsym}   % latex symbols (like \pounds)
 \usepackage{graphicx}   % graphic package for postscript figures
 \usepackage{verbatim}   % for comment environment

\usepackage[all]{xy}   % commutative diagram
\newcommand{\nR}{\mathbb{R}}                     % real number
\newcommand{\nQ}{\mathbb{Q}}                     % complex number

\newcommand{\nP}{\mathbb{P}}                     % projective space

\newcommand{\nA}{\mathbb{A}}                     % affine space

\newcommand{\sO}{\mathscr{O}}                    % structure sheaf
\newcommand{\sI}{\mathscr{I}}                    % ideal sheaf

\newcommand{\mf}[1]{\mathfrak{#1}}

\DeclareMathOperator{\Hom}{Hom}                  % Hom
\DeclareMathOperator{\mld}{mld}                  % mld
\DeclareMathOperator{\Reg}{reg}                  % reg
\DeclareMathOperator{\reg}{reg}                  % reg
\DeclareMathOperator{\Ext}{Ext}                  % Ext
\DeclareMathOperator{\Spec}{Spec}                % Spec
\DeclareMathOperator{\spec}{Spec}                % Spec
\DeclareMathOperator{\Proj}{Proj}                % Proj
\DeclareMathOperator{\proj}{Proj}                % Proj
\DeclareMathOperator{\htt}{ht}                   % ht
\newcommand{\uses}[3]{0\rightarrow {#1} \rightarrow {#2} \rightarrow {#3} \rightarrow 0}

\newtheorem{theorem}{Theorem}[section]
\newtheorem{proposition}[theorem]{Proposition}
\newtheorem{lemma}[theorem]{Lemma}

\newtheorem{conjecture}[theorem]{Conjecture}
\newtheorem*{InvAdj}{Inversion of Adjunction}

\theoremstyle{definition}

\newtheorem{remark}[theorem]{Remark}

\numberwithin{equation}{section}

\begin{document}
\title[A Bound for the Castelnuovo-Mumford Regularity]{A Bound for the Castelnuovo-Mumford Regularity of log canonical varieties}

\author{Wenbo Niu}
\address{Department of Mathematics, Statistics, and Computer
Science, University of Illinois at Chicago, 851 South Morgan Street,
Chicago, IL 60607-7045, USA}

\email{wniu2@uic.edu}

\subjclass[2010]{Primary  14Q20,  13A30; Secondary 14J17, 13D02}

\keywords{Castelnuovo-Mumford regularity, log canonical
singularities}

\date{}
\maketitle
\begin{abstract} In this note, we give a bound for the
Castelnuovo-Mumford regularity of a homogeneous ideal $I$ in terms
of the degrees of its generators. We assume that $I$
defines a local complete intersection  with log
canonical singularities.
\end{abstract}
\section{introduction}

Let $I$ be a homogeneous ideal in a polynomial ring
$R=k[x_0,\ldots,x_n]$ over a field $k$ of characteristic zero.
Consider the minimal free resolution of $R/I$ as a graded
$R$-module,
$$\cdots \rightarrow \oplus_j R(-d_{i,j})\rightarrow \cdots \rightarrow\oplus_j R(-d_{1,j})\rightarrow R\rightarrow R/I\rightarrow 0.$$
The {\em Castelnuovo-Mumford regularity}, or simply {\em
regularity}, of $R/I$ is defined by
$$\reg R/I=\max_{i,j}\{{d_{i,j}-i}\}.$$
The regularity of $I$ is defined by $\reg I=\reg R/I+1$. It measures the complexity of the ideal $I$ and its syzygies. For
more discussion of regularity, see the book of
Eisenbud \cite{Eisenbud:GeoSyz} or the survey of Bayer-Mumford
\cite{Bayer:WhatCanComAlgGeo}.

Suppose that $I$ is generated by homogeneous polynomials of
degrees $d_1\geq d_2\geq\cdots\geq d_t$ and defines a projective
subscheme $X=\proj R/I$ in $\nP^n$ of codimension $r$.  It has been
shown that there is a doubly exponential bound for the regularity of
ideal $I$ in terms of the degrees of its generators. An interesting
question is whether one can find better bounds under some reasonable
conditions on $X$, for instance on its singularities.

If $I_X$ is the saturation of $I$, then $\reg I_X$ is equal to the regularity of the ideal sheaf $\sI_X$ and $\reg \sI_X$ is defined as the minimal number $m$ such
that $H^i(\nP^n,\sI_X(m-i))=0$ for all $i>0$.

The first surprising result was worked out by Bertram, Ein and
Lazarsfeld \cite{BEL} when $X$ is a nonsingular variety. They found a
bound for the regularity of $I_X$ which depends linearly on the degrees of the generators of $I$;
namely
$$\reg R/I_X\leq\sum^r_{i=1} d_i-r.$$
This bound is sharp when $X$ is a complete intersection. Chardin
and Ulrich \cite{CU:Reg} use generic linkage to prove the above
bound in the case when $X$ is a local complete intersection with
rational singularities. Recently, applying multiplier ideal sheaves
and Nadel's vanishing theorem, deFernex and Ein
\cite{Ein:VanishLCPairs} proved that this bound holds in the much
more general situation when the pair $(\nP^n,rX)$ is log canonical.

On the other hand, one can try to bound the regularity of $I$. If $X$ is a local complete
intersection with at most isolated irrational singularities, Chardin
and Ulrich \cite{CU:Reg} gave the following bound:
\begin{equation}\label{b1}\reg R/I\leq\frac{(\dim X+2)!}{2}(\sum^r_{i=1} d_i-r),\end{equation}
which also depends linearly on the degrees of generators. Recently,
in his paper \cite{Fall:Reg}, Fall improved (\ref{b1}) to
$$\reg R/I\leq(\dim X+1)!(\sum^r_{i=1} d_i-r).$$
Starting from this formula and using Bertini's Theorem, Fall also
gave an estimate for the regularity of the defining ideal of any
projective subscheme $X$.

Local complete intersections with rational singularities are
canonical. In light of the work of deFernex and Ein
\cite{Ein:VanishLCPairs} and Chardin and Ulrich \cite{CU:Reg}, it is
natural to ask whether  the bound (\ref{b1}) holds for log canonical
singularities. A scheme of finite type over $k$ is {\em local
complete intersection log canonical} if it is a local complete
intersection with log canonical singularities. In this note, we give
an affirmative answer to this question in the following theorem (as
an easy corollary of Theorem \ref{thm4.2}).

\begin{theorem}\label{Intr:thm}Let $R=k[x_0,\ldots,x_n]$, and let $I=(f_1,\ldots,f_t)$ be
a homogeneous ideal, generated in degrees $d_1\geq d_2\geq \cdots\geq
d_t\geq 1$ of codimension $r$. Assume that $X=\Proj R/I$ is local
complete intersection log canonical and $\dim X\geq 1$. Then
$$\reg R/I\leq \frac{(\dim X+2)!}{2}(\sum_{i=1}^r d_i-r).$$
\end{theorem}

Our main idea relies on the generic linkage method used in
\cite{CU:Reg}. By constructing a generic link $Y$ of $X$, we are
able to reduce the problem to the intersection divisor $Z=Y\cap X$
and then proceed by induction on the dimension. However, for this
approach there are two main problems we need to understand. First we
need to know how to pass singularities from $X$ to $Z$. This is the
hard part of our approach and leads to the study of a flat family of
log canonical singularities. Second we need to control the number
and degrees of the defining equations of $Z$, for which there is a
standard method already.

This note is organized as follows. We explore flat families of
log canonical singularities in section 2. By using Inversion of
Adjunction due to Ein and Musta{\c{t}}{\v{a}} \cite{EM:InvAdjLCI},
we prove the following theorem.

\begin{theorem}\label{2.3} Let $f: Y\rightarrow X$ be a flat morphism
of schemes of finite type over $k$. Assume that $X$ and all fibers
of $f$ are local complete intersection log canonical. Then $Y$ is
local complete intersection log canonical.
\end{theorem}

In section 3, we use the generic residual intersection theory
developed by Huneke and Ulrich \cite{HU88} to pass the log canonical
singularities from $X$ to the intersection divisor $Z$. This is
encoded in the following result.

\begin{proposition}\label{b2}Let $S=\Spec R$ be a regular affine scheme over $k$ and $X\subset
S$ be a subscheme defined by $I=(z_1,\ldots,z_t)$ of codimension $r$.
Construct a generic linkage $J$ of $I$ as follows: let
$M=(U_{ij})_{t\times r}$ be a matrix of variables, $R'=R[U_{ij}]$,
$\alpha=(\alpha_1,\ldots,\alpha_r)=(z_1,\ldots,z_t)\cdot M$ and
$J=[\alpha:IR']$. Let $Z$ be a subscheme of $\spec R'$ defined by
the ideal $J+IR'$. If $X$ is local complete intersection log
canonical, then $Z$ is also local complete intersection log
canonical.
\end{proposition}

In the last section, we use induction to obtain the bound of
regularity. The main idea comes from Chardin and Ulrich
\cite{CU:Reg}.

Some natural questions are pointed out by referees. The first
question is whether it is possible to pass singularities from $X$ to
the link $Y=\spec R'/J$ in Proposition \ref{b2}. The main
difficulty here is that there is no natural morphism from $Y$ to $X$
and therefore we do not know how to pass singularities from $X$ to
$Y$. We may use the morphism from $Y$ to $S$, but we do not know
what kind of fiber it will have. However, it is a really interesting
question and we may propose a conjecture on it.

\begin{conjecture} Assume the hypothesis of proposition \ref{b2}. Set $Y=\spec R'/J$.
If $X$ is local complete intersection log canonical, then $Y$ is
also local complete intersection log canonical.
\end{conjecture}

The second question is, comparing with the work of Chardin, Ulrich
and Fall, can we allow $X$ to have some non-log canonical points and
get a similar  bound to Fall's results? Unfortunately, the method we
use in this note seems unable to solve this problem. Admitting
some non-log canonical points on $X$, we cannot show the
intersection divisor $Z$ has the same property as $X$; this would
be an obstruction to Fall's method. But if we could show $Z$ also
admits non-log canonical points, we may reduce the number of
defining equations of $Z$ by one, and this will lead to Fall's
sharper bound.  Nevertheless, we believe the answer of this question
could be positive, and there will be a better bound under weaker
assumptions. Here we make a conjecture in this direction.

\begin{conjecture} Let $R=k[x_0,\ldots,x_n]$, and $I=(f_1,\ldots,f_t)$ be
an homogeneous ideal, generated in degrees $d_1\geq d_2\geq \cdots\geq
d_t\geq 1$ of codimension $r$. Assume that, except for some isolated
points, $X=\Proj R/I$ is local complete intersection log canonical
and $\dim X\geq 1$. Then
$$\reg R/I\leq (\dim X+1)!(\sum^r_{i=1} d_i-r).$$
\end{conjecture}

The author is grateful to Lawrence Ein and Bernd Ulrich, who offered
many important suggestions and helpful discussions which made this
note possible. The author also thanks referees for their kind
suggestions and patient reading.

\section{Flat Family of Log Canonical Singularities}

In the present section, we study a flat family of local complete
intersection log canonical singularities. We begin by recalling
the definitions of minimal log discrepancy and log canonical
singularities. We mainly follow the approach in \cite[Section
7]{Ein:JetSch}.

Consider a pair $(X,Y)$, where $X$ is a normal, $\nQ$-Gorenstein
variety and $Y$ is a formal finite sum $Y=\sum_i q_i\cdot Y_i$ of
proper closed subschemes $Y_i$ of $X$ with nonnegative rational
coefficients $q_i$.

Let $X'$ be a nonsingular variety which is proper and birational over $X$.
If  $E$ is a prime divisor on $X'$, then $E$ defines a divisor {\em over} $X$.
The image of $E$ on $X$ is called the {\em center} of $E$, denoted by $c_X(E)$.

Given a divisor $E$ over $X$, we choose a proper birational morphism
$\mu:X'\rightarrow X$ with $X'$ nonsingular such that $E$ is a
divisor on $X'$, and such that all the scheme-theoretic inverse
images $\mu^{-1}(Y_i)$ are divisors. The {\em log discrepancy}
$a(E;X,Y)$ is defined such that the coefficient of $E$ in
$K_{X'/X}-\sum_i q_i\cdot \mu^{-1}(Y_i)$ is $a(E;X,Y)-1$. This
number is independent of the choice of $X'$.

Let $W$ be a nonempty closed subset of $X$. The {\em minimal log
discrepancy} of the pair $(X,Y)$ on $W$ is defined by
$$\mld(W;X,Y)=\inf_{c_X(E)\subseteq W}\{a(E;X,Y)\}.$$
If $\mld(p;X,Y)\geq 0$ for a closed point $p\in X$, we say that the
pair $(X,Y)$ is {\em log canonical} at $p$. If $(X,Y)$ is log
canonical at each closed point of $X$, we then say that the pair
$(X,Y)$ is log canonical. If $Y=0$, we just write the pair $(X,Y)$
as $X$.

One important theorem on minimal log discrepancy is Inversion of
Adjunction. It is proved for local complete intersection varieties
by Ein and Musta{\c{t}}{\v{a}}. Since it is used very often in our
proofs, we state it here for the convenience of the reader.

\begin{InvAdj}[{\cite[Theorem1.1]{EM:InvAdjLCI}}] Let $X$ be a
normal, local complete intersection variety, and $Y=\sum_iq_i\cdot
Y_i$, where $q_i\in \nR_+$ and $Y_i\subset X$ are proper closed
subschemes. If $D\subset X$ is a normal effective Cartier divisor
such that $D\nsubseteq \cup_i Y_i$, then for every proper closed
subset $W\subset D$, we have
$$\mld(W;X,D+Y)=\mld(W;D,Y|_D).$$
\end{InvAdj}

The local complete intersection log canonical singularities behave
well in flat families. More specifically, consider a flat family
over a local complete intersection log canonical scheme, where all
fibers are also local complete intersection log canonical. Then we
show that the total space itself is local complete intersection log
canonical.

We start with the case where the flat family has a nonsingular
base.

\begin{proposition} \label{2.4}Let $f: Y\rightarrow X$ be a flat morphism of
schemes of finite type over $k$. Assume that $X$ is nonsingular and
each fiber of $f$ is local complete intersection log canonical. Then
$Y$ is local complete intersection log canonical.
\end{proposition}

\begin{proof} Since  $X$ and all fibers are normal and local complete intersections,
by flatness of $f$, we see that $Y$ is normal and a local complete
intersection (\cite[Section23]{Mat86}). By choosing an irreducible
component of $Y$ and its image, we may assume that $Y$ is a variety
and $f$ is surjective. The question is local. We may assume that
$X=\Spec A$ is affine. Choosing $x\in X$, a closed point defined by
a maximal ideal $m$, $\sO_{X,x}$ is a regular local ring with a
maximal ideal $m_x=(t_1,t_2,\ldots,t_n)$ generated by a regular system
of parameters, where $n=\dim X$. Shrinking $X$ if necessary, we can
extend $t_i$ to $X$ and therefore may assume that
$m=(t_1,\ldots,t_n)\subset A$ generated by a regular sequence. Set
$I_i=(t_1,\ldots,t_i)$. Note that $\sO_{X,x}/(t_1,\ldots,t_i)$ is regular.
By shrinking $X$ if necessary, we may assume further that $A/I_i$ is
regular for each $i=1,\ldots,n$. Let $X_i=\Spec A/I_i$ be subschemes of
$X$ and consider the following fiber product

$$
\begin{CD}
 Y_i @>>> Y\\
 @Vf_iVV @VVfV \\
 X_i @>>> X
\end{CD}
$$
By the flatness of $f_i$ and the assumption that each fiber of $f_i$
is a local complete intersection and normal, we obtain that $Y_i$ is
a local complete intersection and normal for each $i=1,\ldots,n$.

Choose a closed point $y$ on the fiber $Y_x=Y_n$. By the flatness of
$f$, $(t_1,\ldots,t_n)$ is a regular sequence in $\sO_{Y,y}$ and therefore the $t_i$'s define divisors $D_1,\ldots,D_n$ around $y$ in $Y$ such
that
$$Y_i=D_1\cap D_2\cap\cdots\cap D_i,\quad\quad\quad\mbox{for } i=1,\ldots,n.$$
Now by Inversion of Adjunction, we have
\begin{eqnarray*}
\mld(y;Y_n) &=& \mld(y; Y_{n-1},D_n|_{Y_{n-1}})\\
            &=& \mld(y; Y_{n-2},D_n|_{Y_{n-2}}+D_{n-1}|_{Y_{n-2}})\\
            &=& \ldots\\
            &=& \mld(y; Y,D_1+\cdots+D_n).
\end{eqnarray*}
From the assumption that $\mld(y; Y_n)\geq 0$, we get that $\mld(y;
Y)\geq 0$, i.e. $Y$ is log canonical at $y$, which proves the
proposition.
\end{proof}

In the general case in which the flat family has a
singular base, we first resolve the singularities of the base,
and then base change to the situation of nonsingular base. However,
after base change, some extra divisors could be introduced on the
new flat family. This means that we need to consider singularities
of pairs on the new flat family.

\begin{theorem}\label{2.3} Let $f: Y\rightarrow X$ be a flat morphism
of schemes of finite type over $k$. Assume that $X$ and all fibers
of $f$ are local complete intersection log canonical. Then $Y$ is
local complete intersection log canonical.
\end{theorem}

\begin{proof} As in the proof of Proposition \ref{2.4}, we may
assume that $X$ and $Y$ are varieties and $Y$ is normal and a local
complete intersection. We need to show $Y$ is log canonical. Take a
log resolution of $X$, $\mu: \widetilde{X}\rightarrow X$, and
construct the fiber product $\widetilde{Y}=Y\times_X \widetilde{X}$:
$$
\begin{CD}
\widetilde{Y} @>>> Y\\
@VgVV @VVfV \\
\widetilde{X} @>\mu>> X
\end{CD}
$$
By Proposition \ref{2.4}, $\widetilde{Y}$ is local complete
intersection log canonical. Since $X$ is log canonical, we can write
the relative canonical divisor $K_{\widetilde{X}/X}=P-N$, where $P$
and $N$ are effective divisors supported in the exceptional locus of
$\mu$, so that $N=\sum E_i$ where $E_i$ are prime divisors with
simple normal crossings. By base change for relative canonical
diviosrs, we have $K_{\widetilde{Y}/Y}=g^*K_{\widetilde{X}/X}$ and
therefore $K_{\widetilde{Y}/Y}=g^*(P)-g^*(N)$.

Denoting the $F_j$'s as distinct irreducible components of the
$g^*(E_i)$ (note that $g^*(E_i)=g^{-1}(E_i)$ as scheme-theoretical
inverse image of $E_i$), we have $g^*(N)=\sum F_j$. This will be
shown in detail at the beginning of the proof of Lemma
\ref{lemma2.2} below.

Now we let $\pi:Y'\rightarrow \widetilde{Y}$ be a log resolution of
$\widetilde{Y}$ such that
\begin{eqnarray*}
K_{Y'/Y} &=& K_{Y'/\widetilde{Y}}+\pi^*K_{\widetilde{Y}/Y}\\
         &=& A-B + \pi^*g^*P-\sum\pi^*F_i
\end{eqnarray*}
where $A$ is the positive part of $K_{Y'/\widetilde{Y}}$ and $B$ is
the negative part of $K_{Y'/\widetilde{Y}}$ and all prime divisors
in the above formula are simple normal crossings. In order to show
$Y$ is log canonical, it is enough to show that the coefficient of
each prime divisor in $B+\sum \pi^*F_i$ is 1. This is equivalent to
showing that the pair $(\widetilde{Y},g^{-1}N)$ is log canonical,
which is shown in the following Lemma \ref{lemma2.2}.
\end{proof}

\begin{lemma} \label{lemma2.2}Let $f:Y\rightarrow X$ be a flat morphism of varieties
such that $X$ is nonsingular and each fiber of $f$ is local complete
intersection log canonical. Assume that $E_1,\ldots,E_r$ are prime
divisors on $X$ with simple normal crossings. Then the pair $(Y,
\sum_{i=1}^{r} f^{-1}(E_i))$ is log canonical, where $f^{-1}$ means
scheme-theoretical inverse image.
\end{lemma}

\begin{proof} From Proposition \ref{2.4}, $Y$ is local complete
intersection log canonical. Also for each divisor $E_i$, the
scheme-theoretical inverse $f^{-1}(E_i)$ is local complete
intersection log canonical. This implies that
$$\sum_{i=1}^r f^{-1}(E_i)=\sum_{j=1}^s F_j$$
where $F_j$ are distinct irreducible components of the subschemes
$f^{-1}(E_i)$. Note that since $f$ is flat, each $F_j$ only appears
in one $f^{-1}(E_i)$, and if some $F_j$'s are in the same
$f^{-1}(E_i)$ then they are disjoint. Furthermore each $F_j$ is a
Cartier normal divisor on $Y$ with local complete intersection log
canonical singularities. We need to show the pair $(Y, \sum F_j)$ is
log canonical.

We prove this by induction on the dimension of $X$. First
assume that $\dim X=1$. Then $E_1,\ldots,E_r$ are distinct points and
$F_1,\ldots,F_s$ are pairwise disjoint. It is enough to show that for
each $j$, $\mld(F_j;Y,F_j)\geq 0$. Choosing a closed point $p\in
F_j$ of $Y$, by Inversion of Adjunction and the fact that $F_j$ has
log canonical singularities, we have $\mld(p;Y,F_j)=\mld(p;
F_j)\geq0$.

Assume $X$ has any dimension. Since $Y$ is log canonical, it is
enough to show that for each $j$, $\mld(F_j; Y, \sum_{t=1}^s
F_t)\geq 0$. Without loss of generality, we prove this for $F_1$ and
assume that $F_1\subseteq f^{-1}(E_1)$. Choosing any closed point
$p\in F_1$ of $Y$, by Inversion of Adjunction, we have
$$\mld(p; Y, F_1+\sum_{t=2}^sF_t)=\mld(p; F_1 , \sum_{t=2}^sF_t|_{F_1}).$$
For $i=2,\ldots,r$, we set $D_i=E_1\cap E_i$ and note that
$\sum_{t=2}^sF_t|_{F_1}=\sum_{i=2}^r f^{-1}(D_i)$, where $f^{-1}$
means scheme-theoretical inverse image. Now we are in the situation
$$f:F_1\rightarrow E_1,$$
where $E_1$ is nonsingular and $D_2,\ldots,D_r$ are divisors on $E_1$
with simple normal crossings. Then applying induction on $F_1$, we
get that the pair $(F_1, \sum_{t=2}^sF_t|_{F_1})$ is log canonical
and therefore $\mld(p; F_1 , \sum_{t=2}^sF_t|_{F_1})\geq 0$ which
proves the lemma.
\end{proof}

If $f: Y\rightarrow X$ is a surjective smooth morphism, then we
can move singularities freely from $Y$ to $X$. Using the notion of
jet schemes, we have a quick proof for this.

Given any scheme $X$, we can associate the $m$-th jet scheme $X_m$
for any positive integer $m$. The properties of jet schemes are
closely related to the singularities of $X$. We may use jet schemes
to describe local complete intersection log canonical singularities.
The work of Ein and Musta{\c{t}}{\v{a}} shows that if $X$ is a
normal local complete intersection variety, then $X$ has log
canonical singularities if and only if $X_m$ is equidimensional for
every $m$. For more information on jet schemes and their application
to singularities, we refer the reader to \cite{Ein:JetSch}.

\begin{proposition} \label{2.7}Let $f:Y\rightarrow X$ be a smooth surjective morphism of
schemes of finite type over $k$. Then $X$ is local complete
intersection log canonical if and only if $Y$ is local complete
intersection log canonical.
\end{proposition}
\begin{proof} First note that since $f$ is smooth, we have $X$ is normal
and a local complete intersection if and only if $Y$ is normal and a
local complete intersection. Since $f$ is smooth and surjective, for every $m$ we
have an induced morphism between $m$-jet schemes $f_m:Y_m\rightarrow
X_m$, which is smooth and surjective \cite[Remark 2.10]{Ein:JetSch}.
Then $Y_m$ is equidimensional if and only if $X_m$ is
equidimensional. Now by \cite[Theorem 1.3]{EM:InvAdjLCI}, we get the
proposition.
\end{proof}

\begin{remark} \label{rmk2.5} In the proof, if $f$ is smooth but not
surjective, we can only get $f_m:Y_m\rightarrow X_m$ is smooth. Then
equidimensionality of $X_m$ will imply that $Y_m$ is
equidimensional. This means that for a smooth morphism
$f:Y\rightarrow X$, if $X$ is local complete intersection log
canonical then $Y$ is also local complete intersection log canonical
singularities. This is a quick proof for a special case of Theorem
\ref{2.3}.
\end{remark}

\section{Log Canonical Singularities in a Generic Linkage}
In this section, we study the log canonical singularities in a generic
linkage. This could be compared to the work in \cite{CU:Reg}
studying rational singularities in a generic linkage. The $s$-generic
residual intersection theory can be found in \cite{HU88}. Throughout
this section, all rings are assumed to be Noetherian $k$-algebras
and a point on a scheme means a point locally defined by a prime
ideal, not necessarily maximal. All fiber products are over the
field $k$ unless otherwise stated.

Let $S=\Spec R$ be an affine scheme and $X\subset S$ be a
codimension $r$ subscheme defined by an ideal $I=(z_1,\ldots,z_t)$. For
an integer $s\geq 0$, let $M=(U_{ij})_{t\times s}$ be a $t\times s$
matrix of variables and $R'=R[U_{ij}]$ be the polynomial ring over
$R$ obtained by adjoining the variables of $M$. Define $S'=S\times \nA^{t\times
s}=\Spec R'$, which has a natural flat projection $\pi:S'\rightarrow
S$. Let $X'=\pi^{-1}(X)$ be defined by the ideal $IR'$. Construct
an ideal $ \alpha$ in $R'$ generated by $\alpha_1,\ldots,\alpha_s$ as
follows:
$$\alpha=(\alpha_1,\ldots,\alpha_s)=(z_1,\ldots,z_t)\cdot M$$
and set $J=[\alpha:IR']$. The subscheme $Y'$ of $S'$ defined by $J$ is
called an {\em $s$-generic residual intersection} of $X$.

We define $Z$ to be the scheme-theoretical intersection of $X'$ and
$Y'$. Its defining ideal is $I_Z=J+IR'$. We equip $Z$ with a
restricted projection morphism $\pi:Z\rightarrow X$ and call $Z$ an
{\em intersection divisor} of an $s$-generic residual intersection
of $X$.

Note that if $s<r$, then $\alpha$ is generated by a regular sequence
and therefore $J=\alpha$, $Z=X'$. The interesting case is when
$s\geq r$. In particular, when $s=r$, $Y'$ is called a {\em
generic linkage} of $X$. Correspondingly, we call $Z$ an {\em
intersection divisor} of a generic linkage of $X$.

Under the assumption that $X$ is a local complete intersection, the
morphism $\pi:Z\rightarrow X$, and in particular its fibers, can be
understood very well. This offers us an opportunity to pass
singularities from $X$ to $Z$.

We start with a lemma which describes the fibers of $\pi$ when $X$ is
a complete intersection.

\begin{lemma} \label{2.2}Let $S=\Spec R$ be a Gorenstein integral affine scheme and $X$
be a complete intersection subscheme defined by a regular sequence
$I=(z_1,\ldots,z_r)$ in $R$. For $s\geq 0$, let $M=(U_{ij})_{r\times
s}$, $R'=R[U_{ij}]$,
$\alpha=(\alpha_1,\ldots,\alpha_s)=(z_1,\ldots,z_r)\cdot M$ and
$J=[\alpha:IR']$. Assume that $Z$ is defined by $J+IR'$ and consider
the natural morphism $\pi:Z\rightarrow X$. We have
\begin{enumerate}
\item[(1)] If $s<r$, then $Z\cong X\times \nA^{r\times s}$ and $\pi$
is the projection to $X$.
\item[(2)] If $s\geq r$, then $\pi:Z\rightarrow X$ is a flat
morphism and for any point $p\in X$,
$$\pi^{-1}(p)\cong k(p)[U_{ij}]/I_r(M)$$
where $I_r(M)$ is the $r\times r$ minors ideal of $M$.
\item[(3)] In particular, if $s=r$, then $\pi:Z\rightarrow X$ is a flat
morphism such that each fiber is a local complete intersection with
rational singularities.
\end{enumerate}
\end{lemma}
\begin{proof}
(1) is trivial because in this case $J=\alpha$ and $Z$ is defined by
$IR'$ so that $Z=\pi^{-1}(X)\cong X\times \nA^{r\times s}$.

For (2) and (3), picking $\mf{q}\in X\subset S$ and passing to
$R_{\mf{q}}$, we may assume $R$ is local. By  \cite[Example
3.4]{HU88}, $J=(\alpha_1,\ldots,\alpha_s,I_r(M))$. $Z$ is then defined
by $I_Z=J+IR'=(I,I_r(M))$.

Note that
$$R[U_{ij}]/(I,I_r(M))=R/I\otimes_R R[U_{ij}]/I_r(M).$$
This means $\pi:Z\rightarrow X$ can be constructed from the fiber
product
$$
\begin{CD}
Z @>>> \Spec R[U_{ij}]/I_r(M)\\
 @V\pi VV @VV\theta V \\
X @>>> S=\Spec R
\end{CD}
$$
Since $\theta$ is flat, we obtain $\pi$ is flat. The fiber of $\pi$
at $p\in X$ is
\begin{eqnarray*}
F&= &k(p)\otimes_{R/I}R[U_{ij}]/(I,I_r(M))\nonumber\\
&=& k(p)[U_{ij}]/I_r(M) .\nonumber
\end{eqnarray*}
In particular, if $s=r$, we see that $F$ is a local complete
intersection with rational singularities.
\end{proof}

Now we turn to the case where $X$ is a local complete intersection.
\begin{proposition}\label{3.3}Let $S=\Spec R$ be a Gorenstein integral affine scheme and $X$
be a subscheme defined by an ideal $I=(z_1,\ldots,z_t)$ in $R$. For
$s\geq 0$, let $M=(U_{ij})_{t\times s}$, $R'=R[U_{ij}]$,
$\alpha=(\alpha_1,\ldots,\alpha_s)=(z_1,\ldots,z_t)\cdot M$, and
$J=[\alpha:IR']$. Let $Z$ be defined by $J+IR'$ and consider the
natural morphism $\pi:Z\rightarrow X$. Let $\mf{p}\in X$ be a point
of $S$ and assume that $I_{\mf{p}}$ is generated by a regular
sequence of length $r$. Then there is an affine neighborhood of
$\mf{p}$ over which $\pi$ can be factored as follows
$$
\xymatrix{
Z \ar[d]^{\pi}\ar[dr]^{\pi'} &   &  \\
X                         & P\ar[l]^{g} & }
$$
such that $P=X\times\nA^{(t-r)\times s}$ with $g$ the projection to
$X$ and $Z$ can be viewed as an intersection divisor of an
$s$-generic residual intersection of $P$.
\end{proposition}

Note that the above diagram is local. More precisely, there is an
affine neighborhood $U$ of $\mf{p}$ and the morphism $\pi: Z\rightarrow X$
in the above diagram really means the restriction of $\pi$ over $U$, i.e.
$\pi: \pi^{-1}(U)\cap Z\rightarrow U\cap X$.

\begin{proof} By assumption, we may replace $S$ by an affine neighborhood of $\mf{p}$ such that
$I$ is generated by a regular sequence, say $z_1,\ldots,z_r$. Then

\begin{equation}\label{3.1}
\left\{\begin{array}{rcl}
z_{r+1}&= & a_{1,r+1}z_1+a_{2,r+1}z_2+\ldots+a_{r,r+1}z_r\\
z_{r+2}&= & a_{1,r+2}z_1+a_{2,r+2}z_2+\ldots+a_{r,r+2}z_r\\
&\ldots& \\
z_t &= & a_{1,t}z_1+a_{2,t}z_2+\ldots+a_{r,t}z_r
\end{array} \right.
\end{equation}
where $a_{ij}\in R$. Set $A=(a_{ij})_{r\times(t-r)}$. We can write
$(z_{r+1},\ldots,z_t)=(z_1,\ldots,z_r)\cdot A$. Denote
$\displaystyle M={{C} \choose {B}}$, where
$$C= \left(\begin{array}{clcr}
               U_{11} & U_{12}    & \cdots & U_{1s}     \\
               U_{21} & U_{22}    & \cdots & U_{2s}    \\
               \multicolumn{4}{c}{\dotfill} \\
               U_{r1} & U_{r2}    & \cdots & U_{rs}
          \end{array}\right),\ \
B= \left(\begin{array}{clcr}
               U_{r+1,1} & U_{r+1,2}    & \cdots & U_{r+1,s}     \\
               U_{r+2,1} & U_{r+2,2}    & \cdots & U_{r+2,s}    \\
               \multicolumn{4}{c}{\dotfill} \\
               U_{t1} & U_{t2}    & \cdots & U_{ts}
          \end{array}\right).$$
Using the equations in (\ref{3.1}), we can rewrite
$(\alpha_1,\ldots,\alpha_s)=(z_1,\ldots,z_t)\cdot M$ as
$$(\alpha_1,\ldots,\alpha_s)=(z_1,\ldots,z_r)\cdot (A\cdot B+C).$$

Set $N=(V_{lm})_{r\times s}=A\cdot B+C$. Then the ring extension of
$R$ to $R'$ can be obtained by extending twice as follows
$$R\rightarrow R_1=R[U_{ij}|i>r]\rightarrow R'=R_1[V_{pq}].$$
The first extension $R\rightarrow R_1$ gives the morphism $g:\spec
R_1=S\times \nA^{(t-r)\times s}\rightarrow S$. Let
$P=g^{-1}(X)=X\times \nA^{(t-r)\times s}$ defined by $IR_1$ which is
the complete intersection generated by the regular sequence
$(z_1,\ldots,z_r)$ in $R_1$. Restricting $g$ to $P$, we get a
projection $g:P\rightarrow X$. In the second extension,
$R_1\rightarrow R'$, we see that $Z$ can be viewed as an
intersection divisor of an $s$-generic residual intersection of $P$
with morphism $\pi': Z\rightarrow P$.
\end{proof}

Since $Z$ is a generic intersection divisor of $X$, the fibers of the morphism $\pi:Z\rightarrow X$ are local complete intersections with rational singularities and they are log canonical. So the morphism $\pi:Z\rightarrow X$ provides us a flat family of
log canonical singularities, to which results of the previous
section can be applied.

\begin{proposition}\label{3.4}Let $S=\Spec R$ be a regular affine
 scheme and $X$ be a subscheme defined by an ideal $I=(z_1,\ldots,z_t)$ with
codimension $r$ in $S$. Construct a generic linkage $J$ of $I$ as
follows: let $M=(U_{ij})_{t\times r}$, $R'=R[U_{ij}]$,
$\alpha=(\alpha_1,\ldots,\alpha_r)=(z_1,\ldots,z_t)\cdot M$, and
$J=[\alpha:IR']$. Let $Z$ be a subscheme of $\spec R'$ defined by
the ideal $J+IR'$ and consider the natural morphism
$\pi:Z\rightarrow X$. If $X$ is local complete intersection log
canonical, then $Z$ is also local complete intersection log
canonical.
\end{proposition}
\begin{proof}
Choose any point $\mf{p}\in X$. By the assumption, $I_\mf{p}$ is
generated by a regular sequence with length $l\geq r$. By
Proposition \ref{3.3}, there is an affine neighborhood of $\mf{p}$,
over which we can factor $\pi: Z\rightarrow X$ as follows
$$
\xymatrix{
Z \ar[d]^{\pi}\ar[dr]^{\pi'} &   &  \\
X                         & P\ar[l]^{g} & }
$$
such that $P\cong X\times \nA^{(t-l)\times r}$, which is defined by
a regular sequence of length $l$ in $S\times \nA^{(t-l)\times r}$,
and $Z$ is an intersection divisor of a $r$-generic residual
intersection of $P$.

There are two possibilities.

If $l=r$, then by Lemma \ref{2.2}  (3), $\pi':Z\rightarrow P$ is a
flat morphism whose fibers are locally complete intersection log
canonical. Now by Proposition \ref{2.7} and Theorem \ref{2.3} we
obtain that $Z$ is local complete intersection log canonical.

If $l>r$, then by Lemma \ref{2.2} (1), $Z\cong P\times \nA^{l\times
r}$. Using Proposition \ref{2.7}, we have that $Z$ is local complete
intersection log canonical.

\end{proof}

We have passed the singularities from $X$ to $Z$ in above
proposition. As we mentioned in the Introduction, we need to understand
the generators of $Z$. Since $Z$ is defined by $J+IR'$, basically,
we need to know the generators of the generic linkage $J$. The
method we will use here is quite standard in \cite{CU:Reg} and we shall be
brief.

\begin{lemma}\label{lemma3.4} Let $X\subset\nP^n$ be a equidimensional Gorenstein
subscheme with log canonical singularities. Then
$$\Reg \omega_X=\dim X+1,$$
where $\omega_X$ is the canonical sheaf of $X$.
\end{lemma}
\begin{proof}
By assumption, $\omega_X$ is a direct sum of the canonical sheaves of
each irreducible component of $X$. We may assume that $X$ is
irreducible. Since $X$ is log canonical, Kodaira vanishing holds for
$X$ \cite[Corollary 1.3]{kovacs:08}, i.e.
$$H^i(X,\omega_X(k))=0,\quad \mbox{for all } k>0 \mbox{ and } i>0.$$
Note that $H^{\dim X}(X, \omega_X)\neq 0$. Then we see $\reg
\omega_X=\dim X+1$.
\end{proof}

\begin{proposition} \label{4.2}Let $X\subset\nP^n$ be a equidimensional
Gorenstein subscheme with log canonical singularities and
codimension $r$. Assume that $Y\subset\nP^n$ is direct linked with
$X$ by forms of degrees $d_1,\ldots,d_r$. Denote by $J$ the defining
ideal of $Y$ and write $\sigma=\sum_{i=1}^r(d_i-1)$. Then
$J=(J)_{\leq\sigma}$.
\end{proposition}
\begin{proof} Let $I\subset R=k[x_0,\ldots,x_n]$ be the defining ideal of $X$ and $d=\dim
R/I$. Let $b=I\cap J$ be generated by forms in degrees $d_1,\ldots,d_r$
and $\omega$ be the canonical module of $R/I$. If $d=2$, i.e., $X$
is a nonsingular curve, then $(\omega)_{\leq d}=\omega$ by
\cite[Proposition 1.1]{CU:Reg}. If $d>2$, i.e., $\dim X>1$, then
$\reg \omega=\reg\omega_X=d$ by Lemma \ref{lemma3.4} and therefore
we have $(\omega)_{\leq d}=\omega$.

Observe that
$$J/b=\Hom_R(R/I,R/b)=\Ext_R^r(R/I,R)[-d_1-\cdots-d_r]=\omega[d-\sigma].$$
Hence
$(J/b)_{\leq\sigma}=(\omega[d-\sigma])_{\leq\sigma}=(\omega)_{\leq
d}[d-\sigma]=\omega[d-\sigma]$. From the diagram

$$\begin{CD}
0 @>>> (b)_{\leq\sigma} @>>> (J)_{\leq\sigma}@>>>(J/b)_{\leq\sigma}@>>>0\\
  @.@| @VVV @|\\
0 @>>> b @>>> J@>>> J/b@>>>0,
\end{CD}$$
we see $(J)_{\leq\sigma}=J$.
\end{proof}

\section{Bounds for Castelnuovo-Mumford Regularity}

Applying the results we have established, we are able to give a
bound for the Castelnuovo-Mumford regularity of a homogenous ideal which defines a local complete intersection log canonical scheme.
This partially generalizes the work of Chardin and Ulrich
\cite{CU:Reg} and gives a new geometric condition under which a
reasonable bound can be obtained. For the convenience of the reader,
we follow the construction from \cite{CU:Reg} and keep the same
notations.

\begin{proposition}\label{4.3} Let $R=k[x_0,\ldots,x_n]$ and $I\subset R$ be
a homogeneous ideal of codimension $r$ generated by forms
$f_1,\ldots,f_t$ of degrees $d_1\geq d_2\geq \cdots\geq d_t\geq 1$. Let
$$a_{ij}=\sum_{|\mu|=d_j-d_i}U_{ij\mu}x^{\mu},\quad \quad\quad\mbox{for }r+1\leq i\leq t,\  1\leq j\leq r,$$
where $U_{ij\mu}$ are variables. Denote $A=(a_{ij})$,
$K=k(U_{ij\mu})$, $R'=R\otimes_k K$
 and define
$$(\alpha_1,\ldots,\alpha_r)=(f_1,\ldots,f_t){I_{r\times r}\choose A},$$
$J= [(\alpha_1,\ldots,\alpha_r)R': IR']$. Assume that $X=\proj R/I$ is
local complete intersection log canonical. Then $Z=\Proj R'/IR'+J$
is local complete intersection log canonical.
\end{proposition}
\begin{proof} Reduce the question to standard affine covers of
$\nP^n_k$. Without loss of generality, we focus on one affine cover
$U=\Spec R_{(x_0)}$, where $R_{(x_0)}$ means the degree zero part of
the homogeneous localization of $R$ with respect to $x_0$, which is
canonically isomorphic to $k[x_1/x_0,\ldots,x_n/x_0]$.
Set $V=\pi^{-1}(U)$, where $\pi$ is the natural morphism
$\pi:\nP^n_K\rightarrow \nP^n_k$. Note that $V=\Spec R_{(x_0)}'$.
For simplicity, we reset our notations as follows. Replace $R$ by
$R_{(x_0)}$, $R'$ by $R'_{(x_0)}$, $f_i$ by $f_i/x_0^{d_i}$,
and $I$ by $I_{(x_0)}$. Then on the affine open set $U$, $X$ is
generated by $I=(f_1,\ldots,f_t)$ in $R$. We redefine elements of the
matrix $A$ by setting $a_{ij}=\sum
U_{ij\mu}x^{\mu}/x_0^{|\mu|}$. We can see that on $V$, $Z$
is defined by the ideal $J+IR'$, where $J=[\alpha:IR']$ and
 $\alpha=(\alpha_1,\ldots,\alpha_r)$ is defined by the equations in
the assumption. Note that $a_{ij}$'s now become variables over $R$
and therefore $A$ is a matrix of variables over $R$. We restrict
to this affine case in the following proof.

Consider ring extensions $R[a_{ij}]\rightarrow
R[U_{ij\mu}]\rightarrow R'=R\otimes_k K$. The first one is given by
adjoining variables. The second one is the localization of
$R[U_{ij\mu}]$ by the multiplicative set $W=k[U_{ij\mu}]\setminus
\{0\}$. They give morphisms $\phi_1$ and $\phi_2$ respectively:
$$
\begin{CD}
 \Spec R' @>\phi_2>> \Spec R[U_{ij\mu}]@>\phi_1>> \Spec
R[a_{ij}].
\end{CD}
$$
In $R[a_{ij}]$, set $J_1=[\alpha:IR[a_{ij}]]$ and define a subscheme
$Z_1=\spec R[a_{ij}]/(J_1+IR[a_{ij}])$, so that
$Z=(\phi_0\circ\phi_1)^{-1}(Z_1)$. To show $Z$ has the desired
singularities, we just need to show $Z_1$ has the desired
singularities. This is because $\phi_1$ is smooth and it passes
singularities from $Z_1$ to $\phi_1^{-1}(Z_1)$ by Proposition
\ref{2.7}. Our singularities are preserved by localization and so
$\phi_2$ will continue passing singularities from $\phi_1^{-1}(Z_1)$
to $Z$. Hence all we need is to prove the proposition for $Z_1$ in
$\spec R[a_{ij}]$.

To this end, we introduce a new matrix of variables
$B=(b_{lm})_{r\times r}$ and set
$$C={{B}\choose{A\cdot B}}=(c_{uv})_{t\times s},$$
which is also a  matrix of variables over $R$. In the ring
$R[c_{uv}]$, we construct an intersection divisor $Z'$ of $X$ as
follows: let $\alpha'=(\alpha_1',\ldots,\alpha_r')=(f_1,\ldots,f_t)\cdot
C$, $J'=(\alpha':IR[c_{uv}])$ and define $Z'=\Spec
R[c_{uv}]/(J'+IR[c_{uv}])$. Then consider the diagram
$$
\begin{CD}
 \Spec R[a_{ij}] @<q<< \Spec R[a_{ij}]\otimes_k k(b_{lm})\\
 @VVV @VVpV \\
 \Spec R @<<< \Spec R[c_{uv}]
\end{CD}
$$
where $q$ is induced by the base field extension
$R[a_{ij}]\rightarrow R[a_{ij}]\otimes_k k(b_{lm})$, and $p$ is
induced by $R[c_{uv}]\rightarrow R[a_{ij}]\otimes_k k(b_{lm})$,
which is the localization of $R[c_{uv}]$ with respect to the
multiplicative set $k[b_{lm}]\setminus \{0\}$. We note that
$p^{-1}(Z')=q^{-1}(Z_1)=Z_1\otimes_k k(b_{lm})$. By Proposition
\ref{3.4}, $Z'$ is local complete intersection log canonical. Since
$p$ is induced by localization, we obtain that $p^{-1}(Z')$ is also
local complete intersection log canonical. Finally because $q$ is
the base field change of $Z_1$ from $k$ to $k(b_{lm})$, it is easy
to see that $Z_1$ is local complete intersection log canonical if
and only if $q^{-1}(Z_1)=Z_1\otimes_k k(b_{lm})$ is local complete
intersection log canonical. This proves the proposition.
\end{proof}

\begin{theorem}\label{thm4.2} Let $R=k[x_0,\ldots,x_n]$ and $I=(f_1,\ldots,f_t)$ be a
homogeneous ideal, not a complete intersection, generated in degrees
$d_1\geq d_2\geq \cdots\geq d_t\geq 1$ of codimension $r$. Assume that
$X=\Proj R/I$ is local complete intersection log canonical and $\dim
X\geq 1$. Then
$$\reg R/I\leq \frac{(\dim X+2)!}{2}(\sum_{i=1}^r d_i-r-1),$$
unless $R=k[x_0,x_1,x_2]$ and $I=lH$ with $l$ a linear form and $H$
a complete intersection of $3$ forms of degree $d_1-1$, in which case
$\reg R/I=3d_1-5$.
\end{theorem}

\begin{proof}
We construct $R'$, $\alpha=(\alpha_1,\ldots,\alpha_r)$, $J$ and $Z$ as
in Proposition \ref{4.3} and write $\sigma=\sum_{i=1}^{r}(d_i-1)$
and $d=\dim R/I$.

By the assumption that $I$ is not a complete intersection, we may
assume that $d_2\geq 2$. Also we note that if $\sigma=1$, then $\htt
I=1$ and there is a linear form $l$ and a homogeneous ideal $H$
such that $f_i=lh_i$ and $H=(h_1,\ldots,h_t)$, where $h_i$ are all
linear forms, so we get $\reg R/I=\reg R/(l)+\reg R/H=0$. Then we
may assume in the following proof that $\sigma\geq 2$.

We consider the codimension $r$ in two cases.

Case of $r\geq 2$. We proceed by induction on $d$. For $d=2$, we
have $n\geq 3$. Applying \cite[Proposition 2.2]{CU:Reg}, we have
$\reg R/I\leq \frac{(\dim X+2)!}{2}(\sigma-1)$.

Assume that $d\geq 3$. Let $X'=\Proj R'/IR'$ which is local complete
intersection log canonical. Let $(IR')^{top}$ be the unmixed part of
$IR'$; it defines an equidimensional subscheme $X'^{top}$ which
is local complete intersection log canonical and $J$ is directly
linked with $(IR')^{top}$ by $\alpha$. By Proposition \ref{4.2},
$J=(J)_{\sigma}$. Set $Z'=\Proj R'/(IR')^{top}+(J)_{\sigma}$ which
is a Cartier divisor on $X'^{top}$, then in the ring
$R'/(IR')^{top}$, $\overline{J}$ is generated by $d$ forms
$\overline{\beta_1},\ldots,\overline{\beta_d}$ of degrees at most
$\sigma$, which give forms $\beta_1,\ldots,\beta_d$ in $J$ of degrees
at most $\sigma$ such that $Z'=\Proj
R'/(IR')^{top}+(\beta_1,\ldots,\beta_d)$, and therefore we obtain
$Z=\Proj R'/IR'+(\beta_1,\ldots,\beta_d)$. Let
$J'=(\alpha_1,\ldots,\alpha_r,\beta_1,\ldots,\beta_d)$. We have an exact
sequence
$$\uses{R'/IR'\cap J'}{R'/IR'\oplus R'/J'}{R'/IR'+J'}.$$
From this, we get
$$\reg R/I=\reg R'/IR'\leq \max\{\reg (R'/IR'\cap J'),\reg(R'/IR'+J')\}.$$
Since $IR'\cap J'=(\alpha_1,\ldots,\alpha_r)$ is a complete
intersection, $\reg(R'/IR'\cap J')=\sigma$. We just need to bound
$\reg(R'/IR'+J')$. Note that
$IR'+J'=(f_1,\ldots,f_t,\beta_1,\ldots,\beta_d)$ and $\htt (IR'+J')=r+1$.
By assumption of $d_1\geq 2$, we have $\sigma \geq d_{r+1}$.

If $IR'+J'$ is a complete intersection, then some $r+1$ generators
will be a regular sequence. Assume that
$f_{i_1},\ldots,f_{i_p},\beta_{j_1},\ldots,\beta_{j_q}$ are such
generators where $p+q=r+1$. Then
$$\reg R'/IR'+J'=\sum_{\eta=1}^p(\deg f_{i_\eta}-1)+\sum_{\mu=1}^q (\deg \beta_{i_\mu}-1).$$
If $p\leq r$, then we can get $\reg
R'/IR'+J'\leq\sigma+d(\sigma-1)\leq\frac{(d+1)!}{2}(\sigma-1)$.
Otherwise $p=r+1$, then we still have $\reg R'/IR'+J'\leq
\sigma+\sigma-1\leq\frac{(d+1)!}{2}(\sigma-1)$.

If $IR'+J'$ is not a complete intersection, then let
$f_{i_1},\ldots,f_{i_p},\beta_{j_1},\ldots,\beta_{j_q}$ be $r+1$ highest
degree generators. By Proposition \ref{4.3}, $Z=\proj R'/IR'+J'$ is
local complete intersection log canonical, then we use induction for
$IR'+J'$ to get
$$\reg R'/IR'+J'\leq\frac{d!}{2}(\sum_{\eta=1}^p(\deg f_{i_\eta}-1)+\sum_{\mu=1}^q (\deg \beta_{i_\mu}-1)-1).$$
If $p\leq r$, then the left part of the equality is
$\leq\frac{d!}{2}( \sigma +d(\sigma-1)-1)\leq\frac{(d+1)!}{2}(\sigma
-1)$. If $p=r+1$, then the left part is $\leq\frac{d!}{2}( \sigma
+d_{r+1}-1-1)\leq\frac{d!}{2}( \sigma
+\sigma-1-1)\leq\frac{(d+1)!}{2}(\sigma -1)$. Hence we still obtain
$$\reg R'/IR'+J'\leq\frac{(d+1)!}{2}(\sigma -1).$$
This proves the result for the case $r\geq2$.

Case of $r=1$. There is an homogeneous form $l$ and an homogeneous
ideal $H$ such that $f_i=lh_i$, $I=lH$ and $H=(h_1,\ldots,h_t)=[I:l]$.
Since $X$ is a local complete intersection and normal, $\htt H\geq
n$. Also by assumption of $d\geq 2$ we have $n\geq 2$. We consider
the following two cases for $n$.

$n=2$, then $R=k[x_0,x_1,x_2]$, $\htt I=1$. Applying
\cite[Proposition 2.2]{CU:Reg}, we get $\reg R/I\leq 3(\sigma-1)$,
unless $R=k[x_0,x_1,x_2]$, $l$ is a linear form and $H$ a complete
intersection of 3 forms of degree $d_1-1$, in which case $\reg
R/I=3d_1-5$.

$n\geq 3$, then $d=n$. We first note that we have the inequality
$$\deg l+\sum_{i=1}^{n+1}(\deg h_i-1)\leq\frac{(n+1)!}{2}(\sigma-1).$$
If $\htt H=n+1$, then $\dim R/H=0$, and thus we have $\reg
R/H\leq\sum_{i=1}^{n+1}(\deg h_i-1)$, from which we get $\reg
R/I=\reg R/(l)+\reg R/H\leq\frac{(n+1)!}{2}(\sigma-1)$. If $\htt
H=n$ and $H$ is a complete intersection, it is easy to see $\reg
R/I\leq\frac{(n+1)!}{2}(\sigma -1)$. If $\htt H=n$ and $H$ is not a
complete intersection, then by \cite[Proposition 2.1]{CU:Reg}, $\reg
R/H\leq \sum_{i=1}^{n+1}(\deg h_i-1)$. So we still obtain $\reg
R/I\leq \frac{(n+1)!}{2}(\sigma -1)$.
\end{proof}

\begin{remark} It is well known that if $I$ is a complete
intersection, then $\reg R/I\leq\sigma$.  Including the situation
of a complete intersection in the theorem above, we get Theorem 1.1
in the Introduction.
\end{remark}

\end{document}